\input amstex.tex
\documentstyle{amsppt}
\input amsppt1.tex
\nologo
\magnification=1000
\NoBlackBoxes

\define\op#1{\operatorname{#1}}

\topmatter
\author   A. Kazarnovski-Krol          \endauthor
\title Harish-Chandra decomposition for zonal spherical function of
type $A_n$   \endtitle 
\address{
  Department of Mathematics
  Rutgers University
  New Brunswick, NJ 08854, USA}
\abstract{This  paper is devoted to homological treatment of
Harish-Chandra decomposition for zonal spherical functions of type $A_n$.
 }
\endtopmatter

\document

\head 0. Introduction\endhead

Heckman-Opdam system of differential equations 
is holonomic , with regular singularities  and has 
locally $|W|$-dimensional space of  solutions ( cf. corollary 3.9 of [12]),
 where $|W|$
is the cardinality of the Weyl group $W$.
The system is a generalization of radial parts of Laplace-Casimir
operators on symmetric Riemannian spaces of nonpositive curvature and
is isomorphic to Calogero-Sutherland model in the integrable systems.
\smallskip
Harish-Chandra asymptotic solution is a unique solution of the system
with the prescribed asymptotic behavior:
$$F_w(z)=z^{w \lambda+\rho}(1+\ldots)$$
($0<|z_1|<|z_2|< \ldots <|z_{n+1}|$).
Here $w \in W$ are elements of the Weyl group. These solutions provide
a basis in the space of all the solutions in the chamber 
$0 < |z_1|<\ldots< |z_{n+1}|$. Among all the solutions there is a
distinguished one up to the constant multiplier, which admits
continuation to analytic function 
at $z_1=z_2=z_3=\ldots=z_{n+1} \ne 0$.
This solution is referred to as
zonal spherical function. Zonal spherical function is
normalized s.t. it is equal to $1$ at $z_1=z_2=\ldots=z_{n+1}=1$.  
\smallskip
Representation of the  zonal spherical
function as linear combination of elements of the basis (Harish-Chandra
asymptotic solutions) is called Harish-Chandra decomposition.
\smallskip
In ref. [11] we provided an integral representation for the solutions of
Heckman-Opdam system of differential equations in the case of $A_n$. 
We also described contours
for integration $\Delta_w$, integrals over them provide Harish-Chandra
asymptotic solution $F_w(z)$.
In ref. [34] we studied the cycle $\pmb\Delta$ for integration for zonal 
spherical function .
This  paper is devoted to homological treatment of
Harish-Chandra decomposition for zonal spherical functions of type $A_n$.
Namely, we explicitly decompose the distinguished cycle 
$\pmb\Delta$ 
into linear combination of cycles $\Delta_w$ described in [11] and 
check that after 
normalization this turns out to be  the Harish-Chandra decomposition
for zonal spherical function of type $A_n$ (theorem 2.2 and 3.1 below). 
The point of view that linear relations between the
 solutions  reflect
the  linear relations in homology group
 is due to B. Riemann. He also emphasized the importance of the monodromy.
In this case the corresponding homology theory is
 described in [2 , 37].
Harish-Chandra asymptotic solutions correspond to conformal blocks in
conformal field theory ($WA_n$-algebras)  and provide a basis in the
space of 
conformal blocks, zonal spherical function is a particular conformal block, in
the case of $A_2$ see figs.      3a,3b, 3c,3d,3e,3f and 2    below.

\smallskip

\subhead 0.1 Notations \endsubhead

$\alpha_1, \alpha_2, \ldots, \alpha_n$ - simple roots of root system
of type $A_n$

$\Lambda_1,\Lambda_2,\ldots,\Lambda_n$- fundamental weights

$R_{+}$ - set of positive roots

$\delta= {1\over 2} \sum_{\alpha \in R_{+}} \alpha$ -halfsum of
positive roots

$k$- complex parameter ( 'halfmultiplicity' of a root)
$$\rho = {k \over 2} \sum_{\alpha \in R_{+}} \alpha$$

$c(\lambda,k)$- c-function of Harish-Chandra

\head 1. Multivalued form\endhead

Consider the following set of variables:

$z_l, \;\;l=1, \ldots, n+1, \; t_{ij} ,\; i=1, \ldots, j,\;\;
j=1,\ldots, n.$

Variables $z_l$ have meaning of arguments, while variables
$t_{ij}$ are variables of integration. 

It is convenient to organize variables $z_l,\; t_{ij}$ in the
 form of a pattern, cf. fig 1.

\midinsert
$$
\matrix
z_{1}&&z_{2}&&\ldots&&\ldots&&z_{n+1}\\\\
&t_{1,n}&&t_{2,n}&&\ldots&&t_{n,n}&\\\\
&&\ldots&&\ldots&&\ldots&&\\\\
&&&t_{1,2}&&t_{2,2}&&&\\\\
&&&&t_{1,1}&&&&
\endmatrix
$$
\botcaption{Figure 1}
Variables organized in a pattern 
\endcaption
\endinsert

The idea of such an organization is borrowed from Gelfand-Zetlin
patterns [1].

\definition{Definition 1.1} Consider the following multivalued form
$\omega(z,t):$

$$
\align
\omega(z,t):= &\prod_{i=1}^{n+1} z_i^{\lambda_1 +{k n\over 2}} 
 \prod_{i_1 > i_2}{(z_{i_1}-z_{i_2})^{1-2k}}\\
 &\times\prod_{l=1}^{n+1}\prod_{i=1}^n{(z_l-t_{i,n})}^{k-1}\\
 &\times\prod_{j=1}^{n-1} 
 \prod_{i_1=1}^{j+1} \prod_{i=1}^j ( t_{ij} - t_{i_1,j+1} )^{k-1}\\
 &\times\prod_{j=2}^n \prod_{i_1>i_2} {(t_{i_1,j}-t_{i_2,j})^{2-2k} }\\
 &\times\prod_{j=1}^n\prod_{i=1}^{j} {t_{ij}^{\lambda_{n-j+2}-
 \lambda_{n-j+1} -k} } \quad
 {dt_{11} dt_{12} dt_{22} \ldots dt_{nn} }
\endalign
$$

\enddefinition

\remark{Remark 1.2} $k$ is a complex parameter - `halfmultiplicity' of a 
restricted root, cf. Heckman, Opdam [12].
\endremark

In  [11] we proved that integrals over the form $\omega(z,t)$ over
appropriate cycles provide all the solutions to
Heckman-Opdam system of differental equtions and described cycles
$\Delta_w$ for
Harish-Chandra asymptotic solutions ( definition 4.3 and theorem 6.3
of [11]).

\definition{Definition 1.3}
A complex number $z$ can be represented as
$z=r e^{i \alpha}$, where $r, \alpha$ are real numbers, $r \ge 0$. 
$r$ is called
absolute value of $z$, while $\alpha$ is called the phase of $z$. When we say 
that the phase of a complex number $z$ is equal to $0$, we mean
that $\alpha=0$, or the number itself is real and nonnegative. 
\enddefinition
\smallskip
\head 2. The distinguished cycle $\pmb\Delta$ \endhead

Assume that $z_1, z_2, \ldots , z_{n+1}$ are real and 

$$0 < z_1 <z_2 < \ldots <z_{n+1}.$$

\definition{Definition 2.1} Define cycle $\pmb\Delta=\pmb\Delta(z)$ by the following inequalities:

$t_{i,j+1} \le t_{ij} \le t_{i+1, j+1}$ and 

$z_i \le t_{in} \le z_{i+1}$  .

Define form $\omega_{\Delta}(z,t)$ as:

$$
\align
\omega_{\Delta}(z,t):= &\prod_{i=1}^{n+1} z_i^{\lambda_1 +{k n\over 2}} 
 \prod_{i_1 > i_2}{(z_{i_1}-z_{i_2})^{1-2k}}\\
 &\times\prod_{i \le l}{(z_l-t_{i,n})}^{k-1}\prod_{i > l}{(t_{i,n}-z_l)}^{k-1}  \\
 &\times\prod_{j=1}^{n-1} 
 \prod_{i_1 > i_2} ( t_{i_1,j} - t_{i_2,j+1} )^{k-1} 
\prod_{i_2 \ge  i_1} ( t_{i_{2},j+1 } - t_{i_1,j} )^{k-1}\\
 &\times\prod_{j=2}^n \prod_{i_1>i_2} {(t_{i_1,j}-t_{i_2,j})^{2-2k} }\\
 &\times\prod_{j=1}^n \prod_{i=1}^{j} {t_{ij}^{\lambda_{n-j+2}-
 \lambda_{n-j+1} -k} } \quad
 {dt_{11} dt_{12} dt_{22} \ldots dt_{nn} }
\endalign
$$

It is assumed that phases of factors in the formula for
$\omega_{\Delta}(z,t)$
are equal to zero if $k$ and $\lambda_1, \lambda_2, \ldots,
\lambda_{n+1} $ are real. For the homological meaning of the cycle
$\pmb\Delta$ see fig. 2 below and theorem   5.7 of [34].

\enddefinition

\midinsert\vskip 5cm
\includegraphics{zon.eps}
\vskip 2cm \botcaption{Figure 2}
   Zonal spherical function
\endcaption 
\endinsert

In [11] cycles $\Delta_w(z)$ and forms 
$\omega_w(z,t)$ were described,in the case of $A_2$ see 
figs. 3a,3b,3c,3d,3e,3f below. In the case $A_n$ the figures are similar.

\midinsert\vskip 5cm
\includegraphics{1A.ps}
\botcaption{Figure 3a }  
\endcaption 
\endinsert

\vskip 6cm
\midinsert\vskip 6cm
\includegraphics{1BB.ps}
\vskip 2cm \botcaption{Figure  3b}
\endcaption 
\endinsert

\midinsert\vskip 6cm
\includegraphics{1c.eps}
\vskip 5cm \botcaption{Figure  3c}  
\endcaption 
\endinsert

\midinsert\vskip 6cm
\includegraphics{1d.eps}
\vskip 6cm \botcaption{Figure  3d}
\endcaption 
\endinsert

\midinsert\vskip 4cm
\includegraphics{1e.ps}
\vskip 4cm \botcaption{Figure  3e}
\endcaption 
\endinsert

\midinsert\vskip 6cm
\includegraphics{1f.ps}
\vskip 6cm \botcaption{Figure  3f}
\endcaption 
\endinsert

\smallskip
The following theorem explains the relation between 
$\Delta(z)$ and $\Delta_w(z)$.

\proclaim{Theorem 2.2} (Harish-Chandra decomposition)

$$\int_{\Delta(z)} {\omega_{\Delta(z,t)}}= \sum_{w \in S_{n+1}}
b(w,\lambda,k)\int_{\Delta_w(z)}{\omega_w(z)},$$ where

$$
b(w, \lambda,k)={{e^{2 \pi i(\lambda, \delta) } e^{\pi i \; \; l(w)(k-1)}}\over{
(2i)^{n(n+1)\over2} \prod\limits_{\alpha \in  R_{+}}{ \sin( - \pi (w \lambda,
\alpha^{\vee}))}}}
$$

\endproclaim

The theorem is an application of the two following lemmas, also
section 2 of [11] is useful.

\proclaim{Lemma 2.3} (Elementary decomposition)
Let $z_1,z_2,\ldots, z_n$ be real and  $0<z_1<z_2<\ldots<z_n$.
Consider the following integral:

$$\int {t^{a_0 -1}(z_1-t)^{a_1 -1}(z_2-t)^{a_2
-1}\ldots (z_n-t)^{a_n -1}dt}$$

and consider contours $\gamma_1(t), \gamma_2(t), \gamma(t)$, $t\in
[0,1]$
as follows.
$$\gamma(t)= t z_{i-1} +(1-t) z_i$$
$\gamma_1(t)$ is a loop which starts and ends at $z_{i-1}$  and goes 
counterclockwise s.t. the following inequalities are fullfilled for
all $t \in [0,1]$:
 $$z_{i-2}<|\gamma_1(t)| \le z_{i-1}$$
$\gamma_2(t)$ is a loop which starts and ends at $z_{i}$  and goes 
counterclockwise ,s.t. the following inequalities are fullfilled for
all $t \in [0,1]$:
$$z_{i-1}<|\gamma_2(t)| \le z_i$$
as indicated on fig. 4, phases of the factors should be
appropriately chosen.
The following is the specific choice of the phases :
if all  $a_0, a_1,\ldots,a_n$ are real, then the phase of the
integrand along 
$\gamma_1(t),\gamma_2(t), \gamma(t)$  is chosen to be 
zero for small values of $t$.
Then we have the following relation between $\gamma_1, \gamma_2 $ and $\gamma$:
$$\gamma = {{-\gamma_1 }\times{ e^{- \pi i(a_0 + a_1 + \ldots
+a_{i-1})}\over {(2i) \sin \pi(a_0 + a_1 +
\ldots+a_i)}}} +
{{\gamma_2} \times{ e^{- \pi i(a_0 + a_1 + \ldots
+a_i)}\over {(2i) \sin \pi(a_0 + a_1 +\ldots+a_i)}}}
$$
\endproclaim

\midinsert\vskip 0.5cm
\includegraphics{dec.dec.ps}
\vskip 6cm \botcaption{Figure 4}
   Elementary decomposition.
\endcaption 
\endinsert
\proclaim{Lemma 2.4} (Elimination of `wrong' diagrams). Integrals of the
form
$\omega(z,t)$ ($\omega_w(z,t)$), such that contours for integration of $t_{ij},
\; t_{i,j+1}, \; t_{i+1,j+1}$ are shown on fig. 5, provided k is not
an integer,  
 are equal to zero. We suppose that $t_{ij}$ goes from $t_{i,j+1}$ to
$t_{i+1,j+1}$ , $t_{i,j+1}$ goes from $t_{i+1,j+2}$ to  $t_{i+1,j+2}$,
and
$t_{i+1,j+1}$ goes from $t_{i+1,j+2}$ to  $t_{i+1,j+2}$ cf. fig. 5a.  
The same holds true for $t_{i-1,n-1}, t_{i-1,n},
t_{i,n},$ and $z_{i}$ correspondingly, cf. fig. 5b.
\endproclaim

By the 'wrong' diagrams we mean diagrams ,where the two arrows have the
same target, see figs. 5 and 6 of [34] .

\midinsert\vskip 6cm
\includegraphics{dec.elim.1a.ps}
\botcaption{Figure 5a}
     Cycles of this type are homological to zero
\endcaption
\endinsert

\vskip 6cm

\midinsert\vskip 7cm
\includegraphics{dec.elim.1b.ps}
\botcaption{Figure 5b}
     Cycles of this type are homological to zero
\endcaption
\endinsert

\vskip 6cm
\remark{Remark 2.5} Lemma 2.4 is equivalent to quantum Serre's
relations in the form given in    [3], see also [2]. 
\endremark

\head{3. Normalization}\endhead

Let 

$$
\multline
F_w(z) = ({{\prod_{\alpha \in R_{+}} {\Gamma((-w\lambda,
\alpha^{\vee})) \sin ( \pi (-w \lambda, \alpha^{\vee}))\over
\Gamma((-w\lambda,\alpha^{\vee})+k)}}} \\ {\times
e^{-2 \pi i (\lambda, \delta )} e^{-
\pi i (k-1)l(w) }{ \Gamma(k)^{n(n+1)\over2} (2 i)^{n(n+1)\over2}}})^{-1}
\; \int_{\Delta_w(z)} \omega_w(z,t)
\endmultline
$$

Then 

$$F_w(z)= {z^{w\lambda +\rho}}(1+ \ldots)$$  cf. [11] theorem 6.1.

Also, let
$$
F_{\Delta}(z)= {{\Gamma(k) \Gamma(2 k) \ldots \Gamma((n+1)k)}\over
{{\Gamma(k)^{(n+1)(n+2)\over2}}}} \int_{\Delta(z)} \omega_\Delta(z,t)
$$

Then 
$F_{\Delta}(1,1, \ldots,1)=1$, cf. [10] theorem 1.5.
\smallskip

After this normalization theorem 1 reads as usual Harish-Chandra
decomposition  cf. [12,15].

\proclaim{Theorem 3.1} In the above normalization we have:
$$
F_{\Delta}(z)= \sum_{w \in S_{n+1}} c(w \lambda, k) F_w(z),
$$
where $c(w \lambda, k)$ is a $c$-function of Harish-Chandra:

$$c(w \lambda, k) = 
\frac{{\prod_{\alpha \in R_{+}}}\frac{ 
\Gamma((\rho,\alpha^{\vee})+k)}{{\Gamma((\rho,\alpha^{\vee}))}}}
{\prod_{\alpha \in R_{+}}{ 
\Gamma((-w\lambda,\alpha^{\vee})+k)\over{{\Gamma((-w\lambda,\alpha^{\vee}))}}}
}
$$
I.e. $F_{\Delta}(z)$ is identified with zonal spherical function.
\endproclaim

\proclaim{Corollary 3.2} Suppose $ z_1(t), z_2(t), \ldots, z_{n+1}(t)$, $t \in
[0,1]$ are closed loops on a complex plane, i.e. $z_1(0)=z_1(1), 
z_2(0)=z_2(1), \ldots , z_{n+1}(0)=z_{n+1}(1)$, such that $z_i(t) \ne z_j(t)$
for $i \ne j$. Let also $Re(z_i(t)) >0$ for each $i=1,\ldots ,n +1$.
Then the homological class of the cycle $\pmb\Delta$ is preserved under the
monodromy along paths $z_i(t)$.
\endproclaim
\smallskip
\remark{Remark 3.3} In  this approach multiplicative structure of
$c$-function of Harish-Chandra gets a very simple explanation.
Namely:

$$
\multline
c(\lambda,k)=
\frac{
\prod_{1 \le i<j \le n} 
\frac{\Gamma((\rho,e_i-e_j)+k)}{\Gamma((\rho, e_i-e_j))}
      }
 {\prod_{1 \le i<j \le n} 
   \frac{\Gamma((-w \lambda, e_i-e_j) +k)}{\Gamma((-w \lambda,
    e_i-e_j))}}\times
  \frac{
\prod_{1 \le i <n+1} \frac{\Gamma((\rho,e_i-e_{n+1})+k)}
                           {\Gamma((\rho,e_i-e_{n+1}))}}
{\prod_{1 \le i < n+1}\frac{\Gamma((- w \lambda, e_i-e_{n+1})+k)}
                            {\Gamma((-w \lambda, e_i -e_{n+1}))}}
\endmultline
$$
Here $\big\{e_i-e_j|\quad 1 \le i <j \le n+1 \big\}$ are positive roots
of root system of type $A_n$.
Multiplicative properties of $c$-function of Harish-Chandra were
observed by Bhanu-Murti in the case of $SL(n, \Bbb R)$ and in general
case by Gindikin and Karpelevich [17]. $c$-function of Harish-Chandra is
equal to the product of elements of $6 j$-symbols , see fig. 6.
Multilpicative structure of $c$-function of Harish-Chandra amounts to
simple combinatorics related to positive roots , in this case:
$$
\big\{e_i-e_j | 1 \le i <j \le n+1 \big\}=
\big\{e_i-e_j| 1 \le i <j \le n \big \} \bigcup
\big\{e_i-e_{n+1}| 1 \le i <n \big\}.
$$
This combinatorics is both very instructive and restrictive.
\endremark

\midinsert\vskip 5cm
\includegraphics{6.eps}
\vskip 6cm \botcaption{Figure 6 }
$c$-function of Harish-Chandra as a product of elements of $6 j$- symbols. 
\endcaption
\endinsert

\remark{Remark 3.4} We have also checked the monodromy properties of
the cycle $\pmb\Delta$ using quantum group argument, see [34].
\endremark

\remark{Remark 3.5} Harish-Chandra decomposition for zonal spherical
function might be considered as an analogue of
Bernstein-Gelfand-Gelfand  resolution.
\endremark
\subhead{Concluding remark} \endsubhead
We would like to point out once more that the distinguished  cycle
$\pmb\Delta $ appeared in    the classical calculation of Gelfand and Naimark
[16] of zonal spherical function for $SL(n,\Bbb C)$,
originates in the so-called elliptic coordinates and provides a
materialization of the flag manifold.
\subhead Acknowledgments \endsubhead I am grateful to
I.~Gelfand for stimulating discussions concerning  the theory of
spherical functions and the theory of hypergeometric functions, to
S.~Lukyanov for stimulating discussions concerning conformal field
theory, to  V.~Brazhnikov for helpful discussions.

%

\enddocument